\documentclass{article}
\usepackage{graphicx}
\usepackage{amsmath}
\usepackage{amsfonts}
\usepackage{amssymb}

\begin{document}
\[
\]

\begin{center}
\textbf{{A REMARK ON TWO EXTENSIONS OF THE DISC ALGEBRA AND MERGELYAN'S
THEOREM} }
\end{center}

\vspace{4mm}

\begin{center}
V. Nestoridis and I. Papadoperakis
\end{center}

\vspace{8mm}

\begin{center}
Abstract
\end{center}

We investigate the set of uniform limits of polynomials on any closed Jordan
domain with respect to the chordal metric $\chi$ on $\mathbb{C}\cup\{\infty
\}$. We conclude that Mergelyan's Theorem may be extended to the case of
uniform approximation with respect to $\chi$ on closed Jordan domains. Similar
results are obtained if we replace the one point compactification
$\mathbb{C}\cup\{\infty\}$ of $\mathbb{C}$ by another compactification of
$\mathbb{C}$ homeomorphic to the closed unit disc.

\bigskip

AMS Classification Number: Primary 30J99, secondary 46A99, 30E10.

 \vspace{8mm}

\begin{center}
\textbf{{1. Spherical approximation} }
\end{center}

In $[2]$ it has been considered the following generelization $\tilde{A}(D)$ of
the disc algebra. $\tilde{A}(D)$ contains the constant function $f(z)\equiv
\infty$ for all $z\in\bar{D}$, where $\bar{D}$ is the closed unit disc in
$\mathbb{C}$ and all functions $f:D\longrightarrow\mathbb{C}$ holomorphic in
the open unit disc $D$, such that, for every $\zeta\in\partial D$ the limit
$\lim\limits_{\substack{z\rightarrow\zeta\\ z\in D }}f(z)$ exists in
$\mathbb{C}\cup\{\infty\}$. It has also been proved $([2])$ that $\tilde
{A}(D)$ concides with the set of uniform limits with respect to the chordal
metric $\chi$ of all polynomials on $\bar{D}$.

Now let $\Omega$ be an open Jordan domain in $\mathbb{C}$ and $\bar{\Omega}$
its closure. Let $\phi:D\longrightarrow\Omega$ be a Riemann map. According to
a theorem of Caratheodory $([1]),$ $\phi$ extends to a homeomorphism
$\phi:\bar{D}\longrightarrow\bar{\Omega}$. We consider the set of functions
$f\circ\phi^{-1}:\bar{\Omega}\longrightarrow\mathbb{C}\cup\{\infty\}$ for all
$f\in\tilde{A}(D)$. It is easily seen that this set coincides with $\tilde
{A}(\Omega)$, where $\tilde{A}(\Omega)$ is defined as follows. $\tilde
{A}(\Omega)$ contains the function $g(z)\equiv\infty$ on $\bar{\Omega}$ and
the functions $g:\bar{\Omega}\longrightarrow\mathbb{C}\cup\{\infty\}$
continous on $\bar{\Omega}$ such that $g(\Omega)\subset\mathbb{C}$ and
$g_{|\Omega}$ is holomorphic in $\Omega$. \vspace{4mm}

Theorem 1. Under the above assumptions and notation $\tilde{A}(\Omega)$
coincides with the set of uniform limits with respect to the metric $\chi$ of
polynomials on $\bar{\Omega}$.

Proof. Let $P_{n}$ be a sequence of polynomials and $g:\bar{\Omega
}\longrightarrow\mathbb{C}\cup\{\infty\}$ a function such that
$\sup \limits_{z\in\bar{\Omega}}\chi(P_{n}(z),g(z)))\rightarrow0$,
as $n\rightarrow +\infty$. Then
$\sup\limits_{z\in\bar{D}}\chi(P_{n}\circ\phi(z),g\circ
\phi(z)))\rightarrow0$, as $n\rightarrow+\infty$. Since
$P_{n}\circ\phi\in
A(D)$, there exist polynomials $Q_{n}$ so that $\sup\limits_{z\in\bar{D}}%
\chi(P_{n}\circ\phi(z),Q_{n}(z))<\frac{1}{n}$. It follows that $\sup
\limits_{z\in\bar{D}}\chi(Q_{n}(z),g\circ\phi(z))\rightarrow0$, as
$n\rightarrow+\infty$. Thus $g\circ\phi\in\tilde{A}(D)$, which implies that
$g=(g\circ\phi)\circ\phi^{-1}\in\tilde{A}(\Omega)$.

Conversely, let $g\in\tilde{A}(\Omega)$, then $g=f\circ\phi^{-1}$ for some
$f\in\tilde{A}(D)$. Therefore, there exists a sequence of polynomials $P_{n}$
with $\sup\limits_{z\in\bar{D}}\chi(f(z),P_{n}(z)))\rightarrow0$, as
$n\rightarrow+\infty$. It follows that $\sup\limits_{z\in\bar{\Omega}}%
\chi(g(z),P_{n}\circ\phi^{-1}(z)))\rightarrow0$, as $n\rightarrow+\infty$.
Since $P_{n}\circ\phi^{-1}\in A(\bar{\Omega})$, the classical Mergelyan's
Theorem implies that there exist polynomials $Q_{n}$ satisfying $\sup
\limits_{z\in\bar{\Omega}}|P_{n}\circ\phi^{-1}(z)-Q_{n}(z)|<\frac{1}{n}$.
Since for all $a,b\in\mathbb{C}$ we have $\chi(a,b)\leq|a-b|,$ it follows
$\sup\limits_{z\in\bar{\Omega}}\chi(P_{n}\circ\phi^{-1}(z),Q_{n}(z))<\frac
{1}{n}$. The triangle inequality imples $\sup\limits_{z\in\bar{\Omega}}%
\chi(Q_{n}(z),g(z))\rightarrow0$, as $n\rightarrow+\infty$, thus g is the
uniform limit with respect to $\chi$ of the squence of polynomials $Q_{n}$ on
$\bar{\Omega}$. This completes the proof. \vspace{4mm}

\begin{center}
\textbf{{2. Another compactification of $\mathbb{C}$} }
\end{center}

We identify $\mathbb{C}$ with $D$ by the homeomorphism $\mathbb{C}\ni
z\longrightarrow\frac{z}{1+|z|}\in D$. Since $\bar{D}$ is a compactification
of $D$, it induces a compactification $\bar{\mathbb{C}}=\mathbb{C}%
\cup{\mathbb{C}}^{\infty}$, where ${\mathbb{C}}^{\infty}=\{\infty\cdot
e^{i\theta}:\theta\in\lbrack0,2\pi)\}$. The usual Eucledian distance on
$\bar{D}$ iduces a metric $d$ on $\bar{C}$ where $d(z,w)=|\frac{z}%
{1+|z|}-\frac{w}{1+|w|}|$ for $z,w\in\mathbb{C}$, $d(z,\infty\cdot e^{i\theta
})=|\frac{z}{1+|z|}-e^{i\theta}|$ for $z\in\mathbb{C}$, $\theta\in
\lbrack0,2\pi)$ and $d(\infty\cdot e^{i\theta},\infty\cdot e^{i\phi
})=|e^{i\theta}-e^{i\phi}|$ for $\theta,\phi\in\lbrack0,2\pi)$. In $[3]$ it
has been investigated the set of uniform limits with respect to the metric $d$
of the polynomials on $\bar{D}$. This set coincider with the class $\bar
{A}(D)$ defined as follows. $\bar{A}(D)$ contains continuous functions
$f:\bar{D}\longrightarrow\bar{\mathbb{C}}$ of two types. The finite type are
those $f$'s such that $f(D)\subset\mathbb{C}$ and $f_{|D}$ is holomorphic. The
infinite type are those $f$'s such that $f(\bar{D})\subset\mathbb{C}^{\infty}$
and $f(z)=\infty\cdot e^{i\theta(z)}$ where $\theta:\bar{D}\longrightarrow
\mathbb{R}$ is continuous on $\bar{D}$ and harmonic on $D$.

Let $\Omega\subset\mathbb{C}$ be an open Jordan domain and $\bar{\Omega}$ its
closure. We consider $\phi:D\longrightarrow\Omega$ a Riemann map which it is
known that it extends to a homeomorphism $\phi:\bar{D}\longrightarrow
\bar{\Omega}$. We consider the set of functions $f\circ\phi^{-1}:\bar{\Omega
}\longrightarrow\bar{C}$ for all $f\in\bar{A}(D)$. It is easily seen that this
set coincides with the class $\bar{A}(\Omega)$ defined as follows: $\bar
{A}(\Omega)$ contains continuous functions $g:\bar{\Omega}\longrightarrow
\bar{\mathbb{C}}$ of two types. The finite type is those $g$'s with
$g(\Omega)\subset\mathbb{C}$ and $g_{|\Omega}$ holomorphic. The infinite type
is those $g$'s with $g(z)\in{\mathbb{C}}^{\infty}$ for all $z\in\bar{\Omega}$
and $g(z)=\infty\cdot e^{i\theta(z)}$ where $\theta:\bar{\Omega}%
\longrightarrow\mathbb{R}$ is continuous on $\bar{\Omega}$ and harmonic in
$\Omega$. \vspace{4mm}

Theorem 2. Under the above assumptions and notation the class $\bar{A}%
(\Omega)$ coincides with the set of uniform limits with respect to the metric
$d$ of polynomials on $\bar{\Omega}$.

The proof is similar to that of Theorem 1 and is ommitted. \vspace{8mm}

\begin{center}
\textbf{{References} }
\end{center}

[1] Paul Koosis, Introduction to Hp spaces, London Math. Soc. Lecture Note
Series 40, Cambridge University Press, Cambridge, London, N. Y., New Rochelle,
Melbourne, Sydney.

[2] V. Nestoridis, An extension of the disc algebra, Arxiv: 1009.5364

[3] V. Nestoridis and N. Papadatos, Another extension of the disc algebra,
Arxiv: 1012.3674

[4] W. Rudin, Real and complex Analysis, McGraw-Hill, N.Y., St. Louis, San
Fransisco, Toronto, London, Syndney.

\bigskip

\bigskip

V. Nestoridis

Department of Mathematics

University of Athens

Panepistemiopolis

157 84 Athens

Greece

e-mail address: vnestor@math.uoa.gr

\bigskip

I. Papadoperakis

Laboratory of Mathematics

Agricultural University of Athens

118 55 Athens

Greece

e-mail: papadoperakis@aua.gr
\end{document}